




\magnification=1200
\baselineskip=18pt
\def\n{\noindent}
\def\ss{\smallskip}
\def\ms{\medskip}
\def\bs{\bigskip}
\def\lmd{\lambda}
\def\alf{\alpha}
\def\vp{\varepsilon}

\def\cn{{\cal C}_n}
\def\R{\mathop{{\rm I}\kern-.2em{\rm R}}\nolimits }
\def\Rn{\hbox{${\mathop{{\rm I}\kern-.2em{\rm R}}\nolimits}^{\hbox{\rm
n}}$} } 
\def\bbbox{\quad{\vrule height4pt depth0ptwidth4pt}}
\def\dsp{\displaystyle}


\centerline{
\bf Hypercontractivity and comparison of moments of iterated maxima and minima }
\centerline{\bf of independent random variables}
\footnote{}{{\it AMS 1991 Subject Classification:} Primary 60B11, 60E07,
60E15. Secondary 52A21, 60G15}

\bigskip

\centerline{ Pawe{\l} Hitczenko\footnote {*}{Supported in part by an NSF
grant}\footnote{$^\dagger$}{Participant, NSF Workshop in Linear Analysis \&
Probability, Texas A\&M University} }

\centerline{ North Carolina State University } \bigskip

\centerline{Stanis{\l}aw Kwapie\'{n}\footnote {**}{Supported in part 
by Polish KBN grant}}
\centerline{Warsaw University and Texas A\& M University}
\bigskip

\centerline{Wenbo V. Li*{$^\dagger$}}

\centerline{University of Delaware }
\bigskip

\centerline{Gideon Schechtman{$^\dagger$}\footnote {$^\ddagger$}{Supported
in part by
US-Israel Binational Science Foundation}
}
\centerline{Weizmann Institute}
\bigskip

\centerline{Thomas Schlumprecht*\footnote{$^\star$}{Supported in part by a Texas Advanced Research Program grant}
}
\centerline{Texas A\&M University}
\bigskip

\centerline{Joel Zinn*{$^\ddagger$}{$^\star$}}\centerline{Texas A\&M University}

\bigskip
\bigskip
\bigskip

\n {\bf Summary:}
We provide necessary and sufficient conditions for hypercontractivity of
the minima of nonnegative, i.i.d. random variables and of both
the maxima of minima and the minima of maxima for such r.v.'s. It turns out
that the idea of hypercontractivity for minima is closely related to small
ball probabilities and Gaussian correlation inequalities.

\vfill\eject

\beginsection Section 1. Introduction.

In this paper we provide necessary and sufficient conditions for
hypercontractivity of the minima of nonnegative, i.i.d. random variables
of the maxima, the minima and the minima of maxima for such r.v.'s, etc(see the definitions in Section 3). We also give sufficient conditions for hypercontractivity of order statistics. Since questions on ``comparison of moments'' of minima have, apparently, not been considered in the literature, we would like to detail some of our motivation.  The
first motivation for considering such
results is as follows. In a recent paper of de~la~P\~ena, Montgomery-Smith
and Szulga (1994), the authors give a pair of conditions which yield
decoupling (or comparison)
theorems for tail probabilities. One is an $L_p$ decoupling (comparison)
condition for maxima of i.i.d. copies of the variables to be compared. The other is a
hypercontractivity condition, again for maxima of i.i.d. copies of the
``larger" of the two of the r.v.'s. As mentioned the conclusion of their
theorem (Theorem 3.5) is a comparison of the tails of the
r.v.'s of the following type: There exists a constant, $c$, such that
$$P(X>c t)\le c P(Y>t).$$
However, this does not give any information about ``small balls" of the
variables (i.e. the probabilities $P(X\le t)$), since $ P(Y>t)\ge {\dsp {1
/ c}}$ yields a trivial inequality. It is our goal to obtain comparison
results about ``small balls" for norms of some Gaussian vectors. One way to
achieve this is to prove the tail comparisons for minima of independent
copies of the variables with a uniform constant. For example, assume one
has the following holding for all $n$. $$P(\min_{i\le n}X_i>c t)\le c
P(\min_{i\le n}Y_i>t).$$ 
This leads easily to 
$$P(X\le c t)\ge P(Y\le
t),$$ which compares ``small balls". Conversely, if one has an inequality
of the type $$P(X>c t)\le P(Y>t),$$
by raising both sides to a power $n$, one gets an inequality for the
minima. Hence, to use their methods it appears that one needs to consider
i.i.d. maxima
of the variables of interest. Therefore, in our case one is led to consider
i.i.d. maxima of i.i.d. minima.
However, the $L_p/L_q$,   $q>p\ge 1$,
comparison of the
 max min of $X^{\prime}$s to that of $Y^{\prime}$s seems difficult to
handle. In this paper we obtain characterizations of the min and max
hypercontractivity separately, as well as the fact that each of the max min
and min max $L_p$ inequalities are equivalent to the combination of the max
$L_p$ inequality and the min $L_p$ inequality.

The second motivation comes from a well known correlation conjecture for
symmetric, convex sets with respect to a mean zero Gaussian measure $\mu$
on $\Rn$, namely, $$
\mu(A \cap B) \ge \mu (A) \mu (B)
\leqno{(1.1)}
$$
for any symmetric, convex sets $A$ and $B$.

In 1977 L. Pitt (1977) proved that the conjecture holds in $\R^2$. Khatri
(1967) and S\v id\' ak (1967, 1968) proved (1.1) when one of the set is a
symmetric slab (a set of the form $\{x\in \Rn: |(x,u)|\le 1\}$ for some
$u\in \Rn$). For more recent work and references on the correlation
conjecture, see Schechtman, Schlumprecht and Zinn (1995), and Szarek and
Werner (1995).
The Khatri-S\v id\' ak result as a partial solution to the general
correlation conjecture has many applications
in probability and statistics, see Tong (1980). In particular, it is one of
the most important tools discovered recently for the lower bound estimates
of the small ball probabilities, see, for example, Kuelbs, Li and Shao
(1995), and
Talagrand (1994).
On the other hand, the Khatri-S\v id\' ak result only provides the correct
lower bound rate up to a constant at the log level of the small ball
probability. If the correlation conjecture (1.1) holds, then the existence
of the constant of the small ball probability at the log level for the
fractional Brownian motion (cf. Li and Shao (1995)) can be shown. Thus,
from the small ball probability point of view, it is clear that
hypercontractivity for minima, small ball probabilities and the correlation
inequalities are all related, in particular for Gaussian random vectors.
One of our goals in this paper is to expose some of these connections and
in particular to introduce the idea of hypercontractivity for minima to
attack the correlation conjecture and its implication for small ball
probabilities.

The third motivation is related to a weaker form of the correlation
conjecture. To set the notation we let
$\cn$ denote the set of symmetric, convex sets in $\Rn$.

\n Since the correlation conjecture iterates, we consider for $\alf\ge 1$,
\proclaim Conjecture $C_{\alf}$. For any $l,n\ge 1$, and any
$A_1,\cdots,A_l\in\cn $, if $\mu$ is a mean zero, Gaussian measure on
$\Rn$, then $$\mu\bigg(\alf(\bigcap_{i=1}^l
A_i)\bigg)\ge\prod_{i=1}^l\mu(A_i).$$ \par

\n One can restate this (as well as the original conjecture) using Gaussian
vectors in \Rn as follows: for $l,n\ge 1$, and any $A=A_1\times\cdots\times
A_l\subseteq {\R}^{nl}$ let $$\eqalign{\|\cdot\|_A &= \hbox{the norm on
$\R^{nl}$ with the unit ball}\, A,\cr \|\cdot\|_l &=
\hbox{the norm on $\R^n$ with the unit ball}\, A_l.\cr}$$ If
$G,G_1,\cdots,G_l$ are i.i.d. mean zero Gaussian random variables in \Rn,
let $${\cal G}=(G,\cdots,G) \hbox{ and } {\cal H}=(G_1,\cdots,G_l).$$
Then, $C_{\alf}$ can be rewritten as:

\proclaim Restatement of Conjecture $C_{\alf}$. For all $l,n\ge1$, and any
$t>0$, $$\eqalign{\Pr&(\|{\cal G}\|_A\le \alf t) =\Pr({\cal G}\in \alf
t(A_1\times\cdots\times A_l))
\cr &\ge\Pr({\cal H}\in t(A_1\times\cdots\times A_l)) =\Pr(\|{\cal
H}\|_A\le t).}$$\par

\n By taking complements, reversing the inequalities and  raising both sides of the inequality to a power, say $N$,
we get:
$$\Pr(\min_{j\le N}\|{\cal G}^j\|_A> \alf t) \le\Pr(\min_{j\le N}\|{\cal
H}^j\|_A> t).$$\par \bs
\n Again, reversing the inequalities and raising both sides to the power $K$,

$$\Pr(\max_{k\le K}\min_{j\le N}\|{\cal G}^{j,k}\|_A> \alf t)
\le\Pr(\max_{k\le K}\min_{j\le N}\|{\cal H}^{j,k}\|_A> t).$$\par

\n Using the usual formula for $p^{\hbox{th}}$ moments in terms of tail
probabilities we would get:

$$\left\|\max_{k\le K}\min_{j\le N}\|{\cal G}^{j,k}\|_A\right\|_p\le\alf
\left\|\max_{k\le K}\min_{j\le N}\|{\cal
H}^{j,k}\|_A\right\|_p.\leqno(1.2)$$ \bs \n Note that if the conjecture
(1.1) were true then (1.2) would hold with $\alf=1$. Even in the case
$K=N=1$, the best that is known is the above inequality with constant
$\sqrt{2}$. (Of course, if $N=1$, the case
$K=1$ is the same as the case of arbitrary $K$.) To see this first let
$T=:\cup_{l=1}^L T_l=:\cup_{l=1}^L\{(f,l): f\in A_l^{\circ}\}$ where
$A_l^{\circ}$ is the polar of $A_l$. Now define the Gaussian processes
$Y_t$ and $X_t$ for $t\in T_l$ by $Y_{f,l}=f(G)$ and $X_{f,l}=f(G_l)$.
Then, $\sup_{t\in T}Y_t=\max_{l\le L}\|G\|_l$ and $\sup_{t\in T}X_t=
\max_{l\le L}\|G_l\|_l$. We now check the conditions of the
Chevet-Fernique-Sudakov/Tsirelson version of Slepian's inequality (see
also, Marcus-Shepp (1972)). Let $s=(f,p)$ and $t=(g,q)$. If $p=q$,
$(Y_s,Y_t)$ has the same distribution as $(X_s,X_t)$, and hence
$$E|Y_s-Y_t|^2= E|X_s-X_t|^2. $$
If $p\neq q$, then $$E|Y_s-Y_t|^2 \le 2\bigg( EY_s^2+EY_t^2\bigg)=2\bigg(
EX_s^2+EX_t^2\bigg) =2E|X_s-X_t|^2 $$
Therefore, in either case one can use $\sqrt{2}$. Hence, by the version of
the Slepian result mentioned above,
$$E\sup_{t\in T} Y_t\le \sqrt{2} E\sup_{t\in T}X_t.$$

On the other hand the results (mentioned above) of De~La~P\~ena,
Montgomery-Smith and Szulga on decoupling allow one to go from an $L_p$
inequality to a probability inequality if one has one more ingredient,
hypercontractivity. By their results if one can prove that there exists a
constant $\gamma<\infty$ such that for all $K, N$ and symmetric, convex
sets \ms
$$\left\|\max_{k\le K}\min_{j\le N}\|{\cal G}^{j,k}\|_A\right\|_p\le
\gamma\left\|\max_{k\le K}\min_{j\le N}\|{\cal H}^{j,k}\|_A
\right\|_p.\leqno(\hbox{\bf Comparison})$$ and for some $q>p$ and all $K,
N$ and symmetric, convex sets
$$\left\|\max_{k\le K}\min_{j\le N}\|{\cal H}^{j,k}\|_A\right\|_q\le
\gamma\left\|\max_{k\le K}\min_{j\le N}\|{\cal H}^{j,k}\|_A
\right\|_p,\leqno(\hbox{\bf Hyper-contr})$$ then one would obtain for some
$\alf$,
$$\Pr(\min_{j\le N}\|{\cal G}^{j}\|_A>\alf t)\le\alf\Pr(\min_{j\le
N}\|{\cal H}^{j}\|_A>t).$$ This easily implies 
$${\Pr}(\|{\cal G}\|_A>\alf t)\le{\Pr}(\|{\cal H}\|_A>t).$$ Since the
constant outside the probability is now $1$ we can take complements and
reverse the inequality. Now, unraveling the norm and rewriting in terms of
$\mu$ we return to the inequality $C_{\alf}$.

\bigskip
\bigskip

The rest of the paper is organized as follows. Section 2 provides some
basic lemmas and notations. Hypercontractivity for minima and some
equivalent conditions are given in section 3. Section 4 presents
hypercontractivity for maxima in a way suitable for our applications. In
section 5, we combine the results in section 3 and 4 to obtain
hypercontractivity for minmax and maxmin, and comparison results for the
small ball probabilities of possibly different random vectors. We also give a sufficient condition for the comparison of moments of order statistics. In section
6, we apply our results to show that the $\alpha$ symmetric stable random
variables with $0< \alpha \le 2$ are minmax and maxmin hypercontractive, which is strongly connected to the regularity of the $\alf$-stable measure of small balls. In
the case of Gaussian random vectors, we show that the modified correlation
inequality $(C_\alf)$ holds if integrated version of $(C_\alf)$ holds.  In the last
section, we mention some interesting open problems and final remarks.
\bigskip

\beginsection Section 2. Notations and Some Basic Lemmas. \par

For nonnegative i.i.d. r.v.'s $\{ Z_j\}$, let $m_n = m_n(Z)=\min_{j \le n}
Z_j$ and $M_n = M_n(Z)=\max_{j \le n} Z_j$. The $r$-norm of the random
variable $W$ is $$\|W\|_r = (E|W|^r)^{1/r} \quad \hbox{for} \quad r>0 $$

\n We will denote $$x\wedge y = \min\{x,y\},\quad x\vee y = \max\{x,y\}.$$ If
$s<t$ then $s\vee (x\wedge t) = (s\vee x)\wedge t$ and it is denoted by
$s\vee x\wedge t$.

\n Throughout this paper, the numbers $p$ and $q$ are fixed and $0 < p < q
<\infty$.

\ms
\proclaim Lemma 2.1.
Assume $\|W\|_q\le C\|W\|_p$. Then
\item{(a)}
for $\alf=2^{1/(q-p)}C^{q/(q-p)}$, $EW^p\le 2EW^pI_{\{W\le \alf\|W\|_p\}}$,
and \item{(b)} for $0\le\lmd\le 1$, $P(W>\lmd \|W\|_p)\ge\bigl ( {\dsp
{(1-\lmd^p) C^{-p}}}\bigr)^{q/(q-p)}$.

\n {\bf Proof.} (a). Note that
$$
EW^pI_{\{W> \alf\|W\|_p\}}
\le E{\dsp {W^q\over (\alf\|W\|_p)^{q-p}}} \le {\dsp {C^q\|W\|_p^p\over
\alf^{q-p}}} ={1 \over 2}EW^p.
$$
Thus $ EW^pI_{\{W> \alf\|W\|_p\}}\le EW^pI_{\{W\le \alf\|W\|_p\}}$ and
$EW^p \le 2 EW^pI_{\{W\le \alf\|W\|_p\}}$.

\ss
\n (b). The result follows from the Paley-Zygmund inequality $$ EW^p \le
a^p+(EW^q )^{p/q} P^{1-p/q}(W>a) $$ with $a=\lmd \|W\|_p$.
\bbbox
\bs

\proclaim Lemma 2.2.
Let $0<\beta<1$, $0< x,y<1$ then
\itemitem{(a)}
$x\ge \beta^{p/( q-p)}$ and $y^{1/q} \le x^{1/ p}$ imply $(1-x) \le
\beta^{-1} p q^{-1}(1-y)$
\itemitem{(b)} ${p q^{-1}}x\ge y$ implies $(1-x)^{1/ q} \le (1-y)^{1/
p}$.\par \bs
\n {\bf Proof.} (a). We have
$$
1-y \ge 1-x^{q/p}=(1-x)qp^{-1}\eta^{(q-p)/p} \ge (1-x)qp^{-1}x^{(q-p)/p}\ge
(1-x)qp^{-1}\beta $$
where the equality follows from the mean value theorem with $x \le \eta \le 1$.

(b). The conclusion follows from the well known fact $(1-y)^\alf \ge 1-\alf
y$ with $\alf=q/p >1$. \bbbox

\proclaim Lemma 2.3. Fix $0<p\le q<\infty$. Let $\mu$ and $\nu$ be positive
measures on $S$ and
$T$, respectively. If $h\colon\ \Rn \to {\R}_+$ is a measurable function
and $\xi_1,\ldots,
\xi_n$ and $\eta_1,\ldots,\eta_n$ are two sequences of independent r.v.'s
such that
for each $i$ and each $x_1,x_2,\ldots, x_{i-1},x_{i+1},\ldots, x_n$ we have
$$(Eh^q(x_1,\ldots, x_{i-1}, \eta_i,x_{i+1},\ldots, x_n)) ^{1/ q} \le
(Eh^p(x_1,\ldots, x_{i-1},\xi_i, x_{i+1},\ldots, x_n))^{1/p},$$ then
$$(Eh^q(\eta_1,\eta_2,\ldots,
\eta_n))^{1/ q} \le
(Eh^p(\xi_1,\xi_2,\ldots, \xi_n))^{1/ p}.$$\par

\n {\bf Proof.} This follows easily by induction and Minkowski's inequality
$$\left(\int_S\left(\int_T |f(s,t)|^p
\mu(dt)\right)^{q/ p}
\nu(ds)\right)^{1/q} \le \left(\int_T \left(\int_S |f(s,t)|^q
\nu(ds)\right)^{p/q}\mu(dt) \right)^{1/ p}\bbbox .$$ 
\bs

\beginsection Section 3. Hypercontractivity for minima.\par

\proclaim Definition 3.1. We say that a nonnegative random variable $W$ is
$\{p,q\}$-min-$\hfil\break $hypercontractive (with constant $C$), if there
exists $C$ such
that for all $n$, $$\|m_n(W)\|_q\le C\|m_n(W)\|_p.$$ In this case we write
$ W\in \min {\cal H}_{p,q}(C)$.\par

\proclaim Lemma 3.2. If $ W\in \min {\cal H}_{p,q}(C)$ then for each $n$
$$\|m_n( W)\|_p \le K\|m_{2n}( W)\|_p \quad \hbox{with} \quad
K=2^{(2q-p)/p(q-p)} C^{q/( q-p)}.$$\par
\ms
\n {\bf Proof.}
Let $H(t)=H_{p,n}(t)=E(m_n(W)\wedge t)^p$ and note that ${\dsp {H(t)/
t^p}}$ is nonincreasing.
Taking $\alf$ as in the Lemma 2.1,
$$\|m_{2n}\|_p^p\ge Em_{2n}^pI_{m_n\le \alf\|m_n\|_p}=EH(m_n)I_{m_n\le
\alf\|m_n\|_p}.$$
Since ${H(t)/ t^p}$ is nonincreasing,
$$\eqalign{EH(m_n)&I_{m_n\le \alf\|m_n\|_p}=E{\dsp {H(m_n)\over
m_n^p}}m_n^p I_{m_m\le\alf\|m_n\|_p}\ge {H(\alf\|m_n\|_p)\over
(\alf\|m_n\|_p)^p}Em_n^pI_{m_n\le \alf\|m_n\|_p}\cr }. $$ Thus by Lemma
2.1,
$$
\|m_{2n}\|_p^p
\ge {1\over 2}{H(\alf\|m_n\|_p)\over (\alf\|m_n\|_p)^p}\|m_n\|_p^p.$$
Furthermore,
$$H(\alf\|m_n\|_p)\ge Em_n^pI_{m_n\le \alf\|m_n\|_p}\ge{
2^{-1}}\|m_n\|_p^p,$$ which gives the conclusion.\bbbox \ms
The following theorem is a min-analog of a result of De~La~P\~ena,
Montgomery-Smith and
Szulga (1994), proved for maxima, (cf. Theorem 4.4 below).

\proclaim Theorem 3.3. Let $0<p<q$, and let $X,Y$ be r.v.'s\ such that
$X\in \min{\cal H}_{p,q}(C)$ and there exists a constant $B$ such that
$\|m_n(Y)\|_q \le B\|m_n(X)\|_q$ for all $n$. Then $P(X\le \tau t) \le
\delta P(Y\le t)$ for all $t\le t_0= \rho\|X\|_p$ for some constants
$0<\delta<1$, and $\tau ,\rho>0$.

\n {\bf Proof.} By Markov's inequality
$$tP^{1/q} (m_n(Y)>t) \le \|m_n(Y)\|_q \le B\|m_n(X)\|_p.$$ By Lemma 2.1
(b) for each $\lambda$, $0<\lambda<1$, it is $$P^{1/p}(m_n(X) >
\lambda\|m_n(X)\|_p) \ge \left({(1-\lambda^p) C^{-p}}\right)^{q/ p(q-p)} =
D.\leqno (3.1)$$ Hence taking $t=t_n = {B D^{-1}} \|m_n(X)\|_p$ we obtain
$P^{1/ q}(m_n(Y)>t_n) \le P^{1/ p}(m_n(X) > {\lambda B^{-1}}Dt_n)$ which
gives $P^{1/ q} (Y>t_n) \le P^{1/ p}(X>{\lambda B^{-1}D }t_n)$ for all $n$.
By Lemmas 2.1 and  3.2 for each $t_{n+1} \le u\le t_n$, this yields
$$P^{1/ q}(Y>u) \le P^{1/ q}(Y>t_{n+1}) \le P^{1/ p}\left( X>{\lambda D
B^{-1}}t_{n+1}\right) \le P^{1/ p}\left(X> {\lambda D (BK)^{-1}}u\right),$$ where $K$ is as in Lemma 3.2. 
Hence denoting ${\lambda D(BK)^{-1}}$ by $\tau$ we have that $P^{1/ q}(Y>u)
\le P^{1/p}(X>\tau u)$ is satisfied for all $u$ such that $\lim_{n\to
\infty} t_n < u \le t_1 = {B D^{-1}}\|X\|_p$. If $u\le \lim_{n\to \infty}
t_n$, then $P(X>\tau u)=1$ and the above inequality holds true for the
obvious reasons. This inequality and
Lemma 2.2 (a) imply that if $p/q<\beta<1$ and $P(X> \tau u) \ge \beta^{p/(
q-p)}$ then $$P(X\le \tau u) \le
\left(p q^{-1}\beta^{-1}\right) P(Y\le u) \quad \hbox{for}\quad u \le {B
D^{-1}} \|X||_p.$$ Let us observe that by (3.1) for each n
$$P(X>\lambda\|m_{2^n}(X)||_p)\geq D^{p/2^n}$$ and hence by Lemma 3.2
$$P(X>\lambda K^{-n}\|X||_p)\geq D^{p/2^n}.$$ Therefore $P(X>\tau s)\ge
\beta^{p/(q-p)}$ if $s\le \tau^{-1} K^{-n}\lambda \|X\|_p$ and $n$ is such
that $D^{p/2^n}\ge \beta^{p/(q-p)}$.

Thus, for $\tau ={\lambda D( KB)^{-1}}$ we can choose $\delta$ to be any number from the interval
$({pq^{-1}\beta^{-1}},1)$  and then $\rho$ can be taken to be equal
$\min\big\{{B D^{-1}},\tau^{-1}\lambda K^{-n}\big\}$ where $n$ is any
integer such that
$D\ge \beta^{2^n/(q-p)}$. With this choice the assertion of the
theorem is satisfied for all $0\leq t\le t_0=\rho\|X\|_p$.\bbbox
\bigskip
In the above theorem, it is important to note that the constants
$\rho,\tau,\delta$
depend, modulo the given parameters $\beta$ and $\lmd$, only on $p,q,C,B.$
On the other hand, in the following theorem, it is important to note that
each constant appearing in the conditions $(i)-(iv)$ of that theorem depends
only on $p,q$ and the constants from the equivalent conditions, e.g., in particular they do not depend on the
random variable $X$.
This will be useful when considering hypercontractivity of maxima of minima
in section 5.

\proclaim Theorem 3.4. Let $X$ be a nonnegative r.v.\ such that
$\|X\|_q<\infty$ and let $0<p<q$. The following conditions are equivalent
\item{(i)} $X\in \min{\cal H}_{p,q}(C)$ for some $C$, \item{(ii)} there
exist $\vp<1,\ \ \tau , \rho>0$ such that $$P(X\le \tau t) \le \vp P(X\le
t) \quad \hbox{\rm for\ all}\quad t\le t_0=\rho\|X\|_p,$$ \item{(iii)} for
each $\vp>0, \rho>0$, there exists $\tau >0$ such that $$P(X\le \tau t) \le
\vp P(X\le t) \quad \hbox{\rm for \ all}\quad t\le t_0=\rho\|X\|_p,$$
\item{(iv)} there exists $\sigma>0$ such that $$ (E(t\wedge \sigma
X)^q)^{1/q} \le (E(t\wedge X)^p)^{1/p} \quad \hbox{\rm for\ all}\quad t \ge
0.\leqno(3.2)$$

\n {\bf Proof.} (i) $\Rightarrow$ (ii). This implication follows
immediately by Theorem 3.3 applied to $Y=X$.

\n (ii) $\Rightarrow$ (iii). If $\tau ,\vp, t_0$ are as in (ii) then by
induction we obtain for each $n$,
$$P\big(X\le \tau ^n {t_0 t_1^{-1}}t\big) \le \vp^nP \big(X \le {t_0
t_1^{-1}}t\big) \le \vp^n P(X\le t)$$ for all $t\le t_1$, where $t_1$ is
any fixed number such that $t_1 \ge t_0$.

\n (iii) $\Rightarrow$ (iv). For each $t,\sigma$ and $r$, $0<r<1$ we have
$$\eqalign{(E(t\wedge\sigma X)^q)^{1/q} & \le \left(t^qP\left(X > {r
\sigma^{-1}}t\right) + r^qt^qP\left(X\le {rt
\sigma^{-1}}\right)\right)^{1/q}\cr &= t\left(1- (1-r^q) P\left(X \le {rt
\sigma^{-1}}\right)\right)^{1/ q}.}$$ On the other hand $$(E(t\wedge
X)^p)^{1/ p} \ge tP^{1/ p}(X>t) = t(1-P(X\le t))^{1/ p}.$$ Therefore, by
Lemma 2.2~(b) the inequality (3.2) is satisfied if $$P(X\le t) \le {p
q^{-1}}(1-r^q)P \left(X\le {rt \sigma^{-1}}\right).$$ Thus if $\vp = {p
q^{-1}}(1-r^q)$, $\tau $, $t_0$ are as in (iii), then the above inequality
and hence $(3.2)$ is fulfilled for $t\le \tau t_0$ and $\sigma \le r\tau
$.\medskip\n If $t>\tau t_0$ then $(E(t\wedge \sigma X)^q)^{1/ q} \le
\sigma\|X\|_q$ and $(E(t\wedge X^p)^{1/p} \ge (E(\tau t_0 \wedge
X)^p)^{1/p}$ and therefore it is enough to choose $\sigma =
\min\big\{{(E(\tau t_0\wedge X)^p)^{1/p}/ \|X\|_q}, \tau r\big\}$ to have
the inequality $(3.2)$ be satisfied for all $t\geq 0$.

\n (iv) $\Rightarrow$ (i). Applying Lemma 2.3 to
$h(x_1,\cdots,x_n)=x_1\wedge\cdots\wedge x_n,\ \xi_i=X_i$ and
$\eta_i=\sigma X_i$, we get (iv) $\Rightarrow$ (i) with $C=\sigma^{-1}$.
\bbbox

\n {\bf Remark 3.5.} It follows by Theorem 3.4 that
$\{p,q\}$-min-hypercontractivity depends only on the existence of the
$q$-moment and a regularity property of the distribution function at 0,
i.e.\ the following property, which we will call subregularity of $X$ (or,
more precisely, of the distribution of $X$) at 0, $$\lim_{\tau \to 0}
\limsup_{t\to 0} {P(X\le \tau t)\over P(X\le t)}=0.$$

\proclaim Theorem 3.6. Fix $q> 1$. If $\{X_i\}_{i\le n}$ is i.i.d sequence
of nonnegative
r.v.'s with the common
distribution function subregular at 0, and such that $E\|X_1\|^q<\infty$,
then there exists a constant $\sigma$ such that for each $n$, and each
function $h\colon \R_+^n \to \R_+$ which is concave in each variable
separately, we have $$(Eh^q(\sigma X_1,\sigma X_2,\ldots, \sigma
X_n))^{1/q} \le Eh(X_1,X_2,\ldots, X_n).$$

\n {\bf Proof.} By Lemma 2.3 it is enough to prove that there exists
$\sigma>0$, such that for each concave $g\colon\ \R_{+}\to \R_{+},
(Eg^q(\sigma X))^{1/ q} \le Eg(X)$. To see this we first note that by Theorem
3.4, $X$ is $\{q,1\}$-min-hypercontractive, therefore there exists
$\sigma>0$ such that
$(Eh^q_t(\sigma X))^{1/ q}\le Eh_t(X)$ for each $t\ge 0$, where $h_t$ is
given by $h_t(x) = x\wedge t$.

Since for each concave $g\colon\ \R_{+}\to \R_{+}$ there exists a measure
$\mu$ on $\R_{+}$ (the measure $\mu$ is given by the condition $\mu((x,y])
= g'_+(x) - g'_+(y)$ where $g'_+(x)$ is the right derivative of $g$ at $x$)
such that $g = \int_{\R_+} h_t\mu(dt) + g(0)$ the theorem follows by
the Minkowski's inequality. \bbbox

\proclaim Corollary 3.7. If $\{X_i\}$, $h$, $q$ are as in Theorem 3.6 and,
additionally, $h$ is $\alpha$-homogeneous for some $\alpha>0$ (i.e.\ $h(tx)
= t^\alpha h(x))$, then the random variable, $W = h(X_1,X_2,\ldots, X_n)$,
is subregular at $0$.

\n {\bf Proof.} Theorem 3.6 implies that $W$ is
$\{q,1\}$-min-hypercontractive and the result follows by Theorem 3.4.
\bbbox
\bs

\beginsection Section 4. Hypercontractivity of maxima

In this section we treat the case of maxima in a way similar to that of
minima in Section~3. However there are some essential differences which do
not allow us to treat these two cases together.

\proclaim Definition 4.1. We say that a nonnegative r.v.\ $W$ is
$\{p,q\}$-max-hypercontractive if there exists a constant $C$ such that for
all $n$
$$\|M_n( W)\|_q \le C\|M_n( W)\|_p.$$
We will write, $ W\in \max {\cal H}_{p,q}(C)$ in this case.

\proclaim Lemma 4.2. Let $\{X_i\}$ be i.i.d. nonnegative r.v.'s. Then
$${\dsp {nP(X>t)\over 1+nP(X>t)}}\le P(M_n>t)\le nP(X>t).$$

\n {\bf Proof.} The right side is obvious and the left follows by taking
complements and using the inequality, ${\dsp {nu/( 1+nu)}}\le 1-(1-u)^n$.
\bbbox

\proclaim Proposition 4.3.
Let $\{X_i\}$ be i.i.d nonnegative r.v.'s. Then for $a >0$, $$\eqalign{ N
P(M_N\le a) EX^rI_{X>a}\le({\dsp {1-P^N(X\le a)\over
P(X>a)}})EX^rI_{X>a}\le EM_N^r.\cr}\leqno(a)
$$
If $b_N$ satisfies $P(X>b_N)\le{\dsp N^{-1}}\le P(X\ge b_N)$, then $$
\eqalign{{\dsp {1\over 2}}\Bigl(b_N^r+N\int_{b_N}^\infty ru^{r-1}P(X>u)\,
du\Bigr)&\le EM_N^r\le\inf\{a^r+N\int_a^\infty ru^{r-1}P(X>u)\, du\}\cr
&\le b_N^r+N\int_{b_N}^\infty ru^{r-1}P(X>u)\, du } \leqno( b) $$

\ms
\n {\bf Proof.} (a) Let $\tau =\inf\{j\le n: X_j>a\}$. Then
$$\eqalign{E\max_{j\le N}X_j^r&\ge E(X_{\tau};\tau\le
N)=\sum_{j=1}^NE(X_j^rI_{X_j>a}I_{\max_{i<j}X_i\le a})\cr &=\sum_{j=1}^N
P^{j-1}(X\le a)EX^rI_{X>a}=({\dsp {1-P^N(X\le a)\over
P(X>a)}})EX^rI_{X>a}\cr
&\ge N P(M_{N-1}\le a) EX^rI_{X>a}
\cr}$$
(b) To see the right hand inequality, just note that
$$EM_N^r=(\int_0^a+\int_a^\infty) ru^{r-1}P(M_N\ge u)\, du\le
a^r+\int_a^\infty
ru^{r-1}P(M_N\ge u)\, du.$$

\n For the left hand inequality, we again break up the integral as above
and using the defining properties of $b_N$ as well as the monotonicity of
${\dsp {x/( 1+x)}}$ in Lemma 4.2:
$$\eqalign{EM_N^r=(\int_0^{b_N}+\int_{b_N}^\infty) &ru^{r-1}P(M_N\ge u)\,
du\cr &\ge {\dsp {1\over 2}}\int_0^{b_N}ru^{r-1}\, dr + {\dsp {N\over
2}}\int_{b_N}^\infty ru^{r-1}P(X\ge u)\, du.\bbbox}$$

Rychlik (1993) obtained an extension of Theorem 3.5 of De~La~Pe\~na,
Montgomery-Smith and Szulga (cf. Asmar, Montgomery-Smith (1993)). In the
next theorem we put more emphasis on the dependence of the constants.

\proclaim Theorem 4.4. Let $0<p<q$ and let $X\in \max {\cal H}_{p,q}(C)$.
Let $Y$ be a nonnegative r.v. If there exists a constant $D$ such that
$\|M_n(Y)\|_q \le D\|M_n(X)\|_q$ for all $n$, then there are constants
$A,B,\rho$ such that
$$EY^qI(Y>At) \le B^qt^qP(X>t)\quad \hbox{\rm for\ all} \quad t\geq
t_0=\rho\|X\|_p.$$
The constants $B,A,t_0$ can be chosen in the following way:\ if
$0<\lambda<1$, we put \hfil\break
$A = {2^{(q+1)/q }CD \lambda^{-1}}$, $B = A \big({C^p (
1-\lambda^p)}\big)^{1/(q-p)}$ and $t_0 = \lambda\|Y\|_p$.

\n {\bf Proof.} First we note that by the Paley-Zygmund inequality,
$$({\dsp {(1-\lmd^p) C^{-p}}})^{q/(q-p)}\le P(M_n(X)>\lmd\|M_n(X)\|_p)\le n
P(X>\lmd\|M_n(X)\|_p).\leqno(4.1)$$
We next note that by Markov's inequality and the assumptions, for
$\rho=2^{1/q}CD$,
$$P(M_n(Y)\le\rho\|M_n(X)\|_p)\ge P(M_n(Y)\le{{\rho}
(CD)^{-1}}\|M_n(Y)\|_q)\ge{\dsp
{1/ 2}}.
$$
Now, by Proposition 4.3, (a), the assumptions above and (4.1), $$
EY^qI_{Y>\rho\|M_n(X)\|_p}\le 2(CD)^q\|M_n(X)\|_p^q\left({\dsp {C^p
(1-\lmd^p)^{-1}}}\right)^{q/(q-p)}P(X>\lmd\| M_n(X)\|_p). $$
Since $\|M_{2n}(X)\|_p\le 2\|M_n(X)\|_p$, we get by interpolation that $$
EY^qI_{Y>{2\rho t/ \lambda}}\le B^qt^qP(X>t) $$ for
$$
B^q=2^{q+1}\big({CD \lambda^{-1}}\big)^q\left({\dsp {C^p
(1-\lmd^p)^{-1}}}\right)^{q/(q-p)}
$$
as long as $\lambda\|X\|_p\le
2t<\lambda\lim_{n\to\infty}\|M_n(X)\|_p=\lambda\|X\|_\infty.$ If $2t\ge
\lambda\|X\|_\infty$, then since $\rho\ge D$, and
$$\|Y\|_{\infty}=\lim_{n\to\infty}\|M_n(Y)\|_q\le D
\liminf_{n\to\infty}\|M_n(X)\|_q=D\|X\|_\infty,$$ $EY^qI_{Y>{2\rho t/
\lambda}}=0.$ The conclusion follows trivially. \bbbox \bigskip \n {\bf
Remark 4.5.} Theorem 4.4 yields immediately that $$P(Y>At) \le B^qP(X>t)
\quad \hbox{for}\quad t\geq t_0.$$

\proclaim Theorem 4.6. Let $X$ be a nonnegative r.v., $p<q$. The following
conditions are equivalent\medskip
\item{(i)} $X\in \max{\cal H}_{p,q}(C)$ for some $C>0$; \item{(ii)} for each
$\rho>0$ there exists a constant $B$ such that $$EX^qI(X>t) \le
B^qt^qP(X>t) \quad \hbox{\rm for\ all}\quad t\ge t_0=\rho\|X\|_p;$$
\item{(iii)} for each $\rho, \vp>0$ there exists a constant $D>1$
such that
$$D^qP(X>Dt) \le \vp P(X>t) \quad\hbox{\rm for \ all} \quad t\ge
t_0=\rho\|X\|_p;$$
\item{(iv)} there exists a constant $\sigma>0$ such that $$E(t\vee \sigma
X)^q)^{1/ q} \le (E(t\vee X)^p)^{1/ p} \quad \hbox{\rm for \ all} \quad
t\geq 0.$$

\n {\bf Proof.} (i) $\Rightarrow$ (ii). By Theorem 4.4 applied to $Y=X$ we
derive an existence of constants $A, \overline B, \bar t_0$ such that
$$EX^qI(X>t) \le \overline B^qt^qP(X>t)\quad\hbox{for all}\quad t\ge \bar
t_0. $$
Hence for any $t\ge t_0$, where $t_0$ is any number $>0$,
$$\eqalign{EX^qI(X>At) &\le EX^qI(X>A(t\vee \bar t_0)) + EX^qI(t<X\le
A(t\vee \bar t_0))\cr
&\le \overline B^q(t\vee \bar t_0)^q P(X>t\vee \bar t_0) + A^q(t\vee \bar
t_0)^q P(X>t)\cr
&\le\left(1\vee\left({\bar t_0/ t_0}\right)^q\right) (\overline B^q+A^q)
t^qP(X>t).}$$

\n (ii) $\Rightarrow$ (iii). If $t_0,B$ are as in (ii) then for $t\ge t_0$
$$\eqalign{EX^q\ln^+ X/t &= \int^\infty_t E{X^qI(X>s)\over s} ds \le B^q
\int^\infty_t s^{q-1}P(X>s)ds \le q^{-1}{B^q} EX^qI(X>t)\cr &\le
q^{-1}{B^q} t^qP(X>t).}$$
Hence, for any $D>1$ we have
$$(\ln D) D^qt^qP(X>Dt) \le EX^q \ln ^+ X/t \le q^{-1}{B^{2q}} t^q P(X>t)$$
and it is enough to choose $D>1$ such that ${B^{2q}/( q\ln D)}<\vp$.

\n (iii) $\Rightarrow$ (ii). If (iii) holds with $0<\vp <1$, $D>1$ and
$t_0>0$ then by induction\break $P(X>D^nt) \le {\vp^n D^{-nq}} P(X>t)$ for
$t\ge t_0$.
Hence
$$\eqalign{EX^qI(X>t) &= \sum^\infty_{k=0} EX^qI(D^kt<X \le D^{k+1}t) \le
\sum^\infty_{k=0} D^{(k+1)q} t^qP(X>D^kt)\cr &\le \sum^\infty_{k=0}
D^{(k+1)q} t^q \cdot {\vp^k D^{-kq}} P(X>t) = {D^q( 1-\vp)^{-1}}
t^qP(X>t).}$$

\n (ii) and (iii) $\Rightarrow$ (iv). Assume that (ii) and (iii) are
fulfilled with constants $B,D,\vp,t_0$.

By (ii) we obtain for $t\ge t_0\sigma$, where at the moment $\sigma$ is any
number $<1$,
$$\eqalign{(E(t\vee \sigma X)^q)^{1/q} &\le t \biggl(P\left(X<{t
\sigma^{-1}}\right) + {\sigma^q t^{-q}} E\left(X^qI\left(X\ge
{t\sigma^{-1}}\right)\right) \biggr)^{1/ q}\cr &\le t\left(1+(B^q-1)
P\left(X \ge{t\sigma^{-1}}\right)\right)^{1/ q}}.$$ On the other hand for
any $R>1$
$$
(E(t\vee X)^p)^{1/p} \ge t(1+(R^p-1) P(X>Rt))^{1/ p}. $$ Hence by Lemma 2.2
(b) the inequality in (iv) holds if $${p q^{-1}} {(B^q-1)( R^p-1)^{-1}}
P\left(X\ge {t\sigma^{-1}}\right) \le P(X>Rt).$$ Therefore if we choose $R$
so that
$$
{p q^{-1}} {(B^q-1)( R^p-1)^{-1}}< {D^q/\vp}, \quad {R\|X\|_p/ 2} \ge t_0
\quad \hbox{ and } \quad
\sigma<{(RD)^{-1}},
$$
then the inequality in (iv) is satisfied for all $t\ge {\|X\|_p/2}$.

If $t < {\|X\|_p/2}$ then $(E(t\vee \sigma X)^q)^{1/q} \le t+\sigma
\|X\|_q$ and $(E(t\vee X)^p)^{1/p} \ge \|X\|_p$ and therefore if
additionally $\sigma < {\|X\|_p( 2\|X\|_q)^{-1}}$ then the inequality in
(iv) is satisfied for all $t\geq 0$.

\n (iv) $\Rightarrow$ (i). This implication is proved in the same way as
the one in Theorem 3.4. It is enough to replace $\wedge$ by $\vee$
everywhere. \bbbox

As in Theorem 3.4 for fixed $p,q$ the constants in the conditions
$(i)-(iv)$ depend only on themselves.

\n {\bf Remark 4.7.}\item{(i)} The equivalence of (ii) and (iii) in Theorem 4.6 can be deduced from more general results (cf. Bingham, Goldie and Teugels).
\item{(ii)} It follows from Theorem 4.4 that if $X$ is
$\{p,q\}$ max-hypercontractive then for some $\vp>0$ and all $r<q+\vp$, $X$
is also $\{r,q+\vp\}$-max-hypercontractive. \item{(iii)}
The property of $\{p,q\}$-max-hypercontractivity is equivalent to
$$\limsup_{D\to \infty} \limsup_{t\to \infty} {D^qP(X>Dt)\over P(X>t)}=0$$
which we will call $q$-subregularity at
$+\infty$.

\proclaim Theorem 4.8. If a nonnegative r.v.\ $X$ is $q$-subregular at
$+\infty$ then there exists a constant $\sigma$ such that for each $n$ and
each $X_i$, $i=1,\ldots, n$ independent copies of $X$, $$(Eh^q(\sigma
X_1,\sigma X_2,\ldots, \sigma X_n))^{1/q} \le Eh(X_1,X_2,\ldots, X_n)$$
for each function $h\colon \R_+^n\to \R_+$ which is in each variable
separately nondecreasing and convex, and $\lim\limits_{x_i\to +\infty}
\big(x_i {\partial h\over \partial x
_i} (x) - h(x)\big)\le 0$.

\n {\bf Proof.} The proof is the same as in the case of Theorem 3.6, except
that we have to replace everywhere $\wedge$ by $\vee$ and that the measure
$\mu$ is given by $\mu((x,y]) = g'_+(y) -g'_+(x)$ and then $$g(x) =
\int_{\R_+} h_t(x) \mu(dt) + \lim_{t\to\infty} (g(t) - tg'(t)). \bbbox $$

In analogy to Corollary 3.7 we obtain

\proclaim Corollary 4.9. If $\{X_i\}$, $h$ are as in Theorem 4.8, $q>1$ and
in an
addition $h$ is $\alpha$-homogeneous for some $\alpha>0$, then the random
variable $W = h(X_1,X_2,\ldots, X_n)$ is $q$-subregular at $+\infty$.\par
\bs
\beginsection Section 5. Hypercontractivity of minmax and maxmin.

In this section we will impose on $X$ both the condition of subregularity
at 0 and that of $q$-subregularity at $+\infty$.

\proclaim Lemma 5.1. Assume that $X$ is a nonnegative random variable which
is subregular at 0 and $q$-subregular at $+\infty$. Then for each $0<p<q$
there exists a constant $\sigma>0$ such that for each $0<s<t<\infty$
$$(E(s\vee \sigma X\wedge t)^q)^{1/q} \le (E(s\vee X\wedge
t)^p)^{1/p}.\leqno (5.1)$$

\n {\bf Proof.} Let $R>1$ be any fixed number, and let $r=R^{-1}$. Let
$t_0$ be any positive number, and let $\tau$ be such that the inequality in
Theorem 3.4 (iii) holds for $\vp = {pq^{-1}}(1-r^q)$ for all $t\le t_0$.
Then let $\alpha = 2^{-1}(E(\tau\wedge ({X/t_0}))^p)^{1/p}$. The constant
$B$ is such that the inequality in
Theorem 4.6 (ii) is true for all $t\ge t_0$ and let $D$ be such that the
inequality in Theorem 4.6 (iii) is satisfied for $$\vp = {q p^{-1}}(R^p-1)
\left({(B^q-1)^{-1}} \wedge { (R^q-1)^{-1}}\right) \quad \hbox{for}\quad
t\ge t_0.$$
We will show that for $\sigma=\min\big\{{\alpha t_0 /\|X\|_q}, {\alpha/ D},
{r/ D}, r\tau \big\}$ the inequality (5.1) holds true for each
$0<s<t<\infty$. Consider the following five cases. \ms \n Case 1. $s \le
\alpha t_0, t\ge \tau t_0$. We have $$ (E(s\vee \sigma X\wedge t)^q)^{1/q}
\le \alpha t_0 +\sigma\|X\|_q $$ and
$$
(E(s\vee X\wedge t)^p)^{1/ p} \ge (E(\tau t_0 \wedge X)^p)^{1/ p}. $$ Since
$\sigma < {\alpha t_0/\|X\|_q}$ the inequality holds by the choice of
$\alpha$. \ms

\n Case 2. $t\le \tau t_0$, $rt>s$. We have $$\eqalign{(E(s\vee \sigma
X\wedge t)^q)^{1/ q} &\le \left(t^qP\left( X>{r \sigma^{-1}}t\right) +
r^qt^qP\left( X \le {rt \sigma^{-1}}\right) \right)^{1/q}\cr &\le t \left(1
+ (r^q-1)P\left( X\le {rt \sigma^{-1}}\right)\right)^{1/ q} } $$ and
$$
(E(s\vee X\wedge t)^p)^{1/ p} \ge t(1-P(X\le t))^{1/p}. $$ Therefore by
Lemma 2.2 (b) the inequality $(5.1)$ holds if $$ {p q^{-1}}(1-r^q) P\left(X
\le {rt \sigma^{-1}}\right) \ge P(X\le t) $$ which is true by the choice of
$\tau$ since $\sigma <r\tau $, $t\le \tau t_0$. \ms

\n Case 3. $t\le \tau t_0$, $rt\le s$. We have $$ \eqalign{(E(s\vee \sigma
X\wedge t)^q)^{1/q} &\le \left(s^qP\left( X \le {s\sigma^{-1}}\right) +
t^qP(\left(X>{s \sigma^{-1}}\right) \right)^{1/ q}\cr &= t \left(1 +
\left(\left({s/t}\right)^q-1\right) P \left(X\le {s
\sigma^{-1}}\right)\right)^{1/q} }
$$
and
$$
\eqalign{ (E(s\vee X\wedge t)^p)^{1/p} &\ge (s^pP(X\le t) +
t^pP(X>t))^{1/p} \cr &= t\left(1 + \left(\left({s/ t}\right)^p-1\right)
P(X\le t)\right)^{1/ p}.} $$
Therefore by Lemma 2.2 (b) to have $(5.1)$ it is enough to show $${p
q^{-1}} {(1-({s/ t})^q)( 1-({s/ t})^p)^{-1}} P\left(X\le {s
\sigma^{-1}}\right) \ge P(X\le t).$$ Since the function ${(1-x^q)/(
1-x^p)}$ is increasing on $\R_+$ and $s/t\ge r$ it is enough to prove that
${p q^{-1}} {(1-r^q)( 1-r^p)^{-1}} P\big(X\le {rt \sigma^{-1}}\big) \ge
P(X\le t)$ which was proved in the preceding case, because $ 1-r^p<1$. \ms
\n Case 4. $s>\alpha t_0$, $t>Rs$. We have $$ \eqalign{(E(s\vee \sigma
X\wedge t)^q)^{1/q} &\le \left(s^qP\left(X\le {s \sigma^{-1}}\right) +
\sigma^q EX^qI\left(X>{s \sigma^{-1}}\right)\right)^{1/q}\cr &\le s\left(1+
(B^q-1) P\left(X>{s\sigma^{-1}}\right)
\right)^{1/q},} $$
which follows by the choice of $B$, since ${s\sigma^{-1}} > {\alpha t_0
\sigma^{-1}}\ge t_0$, and
$$
\eqalign{(E(s\vee X\wedge t)^p)^{1/p}&\ge (s^pP(X\le Rs)+(Rs)^p
P(X>Rs))^{1/p}\cr & = s(1+(R^p-1) P(X>Rs))^{1/ p}.} $$
By Lemma 2.2 (b) it is enough that
$${p q^{-1}} {(B^q-1)( R^p-1)^{-1}} P\left(X>{s \sigma^{-1}}\right) \le
P(X>Rs)$$ since
$\sigma \le {(\alpha\wedge r)/ D}$ it is enough to show $${P\big(X>D\big({
\alpha^{-1}}\vee R\big)s\big)\over P(X>Rs)} \le {q p^{-1}} {(R^p-1)(
B^q-1)^{-1}}=\vp.$$ But, then by the choice of $D$ we have
$${P\big(X>D\big({\alpha^{-1}}\vee R\big)s\big)\over P(X>Rs)} \le
{P\big(X>D\big({\alpha^{-1}}\vee R\big)s\big)\over
P\big(X>\big({\alpha^{-1}}\vee R\big)s\big)}
{P\big(X>\big({\alpha^{-1}}\vee R\big)s\big) \over P(X>Rs)}\le {\vp\over
D^q}<\vp,$$
because $\big({\alpha^{-1}}\vee R\big)s \ge t_0$. \ms \n Case 5. $s>\alpha
t_0, t\le Rs$. We have $$ \eqalignno{(E(s\vee \sigma X\wedge t)^q)^{1/ q}
&\le \left(s^q P\left(X\le {s \sigma^{-1}}\right)
+t^qP\left(X>{s\sigma^{-1}}\right) \right)^{1/q} \cr &= s\left(\!1 +\left(
\left({t/s}\right)^q\!-1\right) P\left(X >
{s\sigma^{-1}}\right)\right)^{1/q} } $$
and
$$
\eqalignno{ (E(s\vee X\wedge t)^p)^{1/p} &\ge (s^pP(X\le
t)+t^pP(X>t))^{1/p}\cr & = s\left(1+ \left(\left({t/ s}\right)^p-1\right)
P(X>t)\right)^{1/ p}.} $$
By Lemma 2.2 (b) it is enough to prove
$$
{pq^{-1}} {(({t/s})^q-1)(({t/ s})^p-1)^{-1}} P\left(X>{s
\sigma^{-1}}\right) \le P(X\ge Rs).
$$
Since ${t/ s} \le R$ it suffices to show that $$ {pq^{-1}}
{(R^q-1)(R^p-1)^{-1}} P\big(X >{s \sigma^{-1}}\big) \le P(X\ge Rs) $$ which
is shown in the same way as in the preceding. \bbbox

Taking into account the remarks before the proof of Theorem 3.4 and after
Theorem 4.6 we check easily that given $p,q$ the constant $\sigma$ depends
only on the min and max hypercontractivity constants of $X$. More exactly
if we define $\min\max {\cal H}_{p,q}(C)$ as the intersection of $ \min
{\cal H}_{p,q}(C)$ and $\max {\cal H}_{p,q}(C)$ then for a random variable
$X\in \min\max {\cal H}_{p,q}(C)$ the constant $\sigma$ depends only on
$p,q$ and $C$.

\proclaim Corollary 5.2. If $X, p,q$ are as in Lemma 5.1, then there exists
a constant $C$ such that if $(X_i)$, $i=1,\ldots, n$ is a sequence of
independent copies of $X$ and $X^{k,n}$ denotes the $k$-th order statistics
of the sequence $(X_i)$, $i=1,\ldots, n$ then $\|X^{k,n}\|_q \le
C\|X^{k,n}\|_p$ and $X^{k,n}$ is $q$-subregular at $+\infty$ and subregular
at 0.

\n {\bf Proof.} The statistic $X^{k,n}$ can be written as
$h(X_1,X_2,\ldots, X_n)$ where for each $i$ and each fixed $x_i,\ldots,
x_{i-1},x_{i+1},\ldots, x_n$ the function $f(x_i) = h(x_1,x_2,\ldots,
x_{i-1}, x_i, x_{i+1},\ldots, x_n) = (s\vee x_i\wedge t)$ for some $0<s<t$,
and all $x_i\in \R_+$. And therefore the first part of the corollary
follows by the observation. The second part is obtained easily because we
have that $X^{k,n}$ is $\{q,p\}$-max and min hypercontractive. \bbbox

The preceding corollary can be considerably generalized. At first let us
define a class ${\cal F}$ of
functions $g\colon \R_+ \to \R_+$ which can by written as $g(x) =
\int_\Delta h_{s,t}(x) \mu(ds,dt)$ for some positive measure $\mu$ on
$\Delta = \{(s,t) \in \R_+\times \R_+\colon \ s\le t\}$ and where $h_{s,t}$
are functions defined by $h_{s,t}(x) = s\vee x\wedge t$. It is
possible to give intrinsic description of functions in ${\cal F}$. Instead
let us observe that if $f$ is twice continuously differentiable on $\R_+$
then $f\in {\cal F}$ if and only if for each $x\in \R_+$, $0\le xf'(x) \le
f(x)$ and $f(0) \ge \int_{\R_+} x(f''(x)\vee 0)dx$. In this case the
measure $\mu$ is given by the following condition:\ for measurable $h\colon
\ \Delta \to \R_+$
$$\int_\Delta h(s,t) \mu(ds,dt) = \int_{\R_+}\left(\sum_{(s,t)\in I(y)}
h(s,t)\right) dy$$
where $I(y)$ is the countable family of open, disjoint intervals with the
union equal\break $\{s\in \R_+\colon \ f'(s) >y\}$. It is not difficult to
prove that we have the representation $f(x) = \int_\Delta h_{s,t}(x)
\mu(ds,dt) + c$ where $c=f(0) - \int_{\R_+} x(f''(x)\vee 0)dx$.

\proclaim Theorem 5.3. If $X$ is a $q$-subregular at $+\infty$ and
subregular at 0 then there exists a constant $\sigma>0$ such that for each
$n$ and each $h\colon \R_+^n\to \R_+$ which in each variable, separately is
in class ${\cal F}$ we have $$(Eh^q(\sigma X_1,\sigma X_2,\ldots, \sigma
X_n))^{1/ q} \le Eh(X_1,X_2,\ldots, X_n).$$
Moreover if $h$ is $\alpha$-homogeneous for some $\alpha>0$ then
$h(X_1,\ldots, X_n)$ is $q$-subregular at $+\infty$ and subregular at 0.

\n {\bf Proof.} The proof follows the same pattern as proofs of Theorems
3.6, 4.8 and Corollaries 3.7, 4.9, and is based on Lemma 5.1 \bbbox

Applying comparison results of Theorems 3.3 and 4.4 we obtain easily

\proclaim Theorem 5.4. Let $X,Y$ be nonnegative r.v.'s such that $X$ is
$\{p,q\}$-min and max hypercontractive. If there exists a constant $B$ such
that $\|m_n(Y)\|_q \le B\|m_n(X)\|_q$ and $\|M_n(Y)\|_q \le B\|M_n(X)\|_q$
for all $n$, then there exists a constant $D$ such that $P(Y\le t)\ge
P(DX\le t)$ for all $t\in \R_+$.

Finally we have

\proclaim Theorem 5.5. If $X\in \min{\cal H}_{p,q}\cap\max{\cal H}_{p,q}$
then, there exists a constant $D$ such that for all $l$ and all
$n_1,k_1,n_2,k_2,\dots n_l,k_l$, $$
\|M_{n_1}m_{k_1}M_{n_2}m_{k_2}\dots M_{n_l}m_{k_l}(X)\|_q\le
D\|M_{n_1}m_{k_1}M_{n_2}m_{k_2}\dots M_{n_l}m_{k_l}(X)\|_p . $$

\n {\bf Proof.} If $X$ is both min-hypercontractive and
max-hypercontractive then, by Remarks 3.5and 4.7(ii) it is subregular at 0
and q -- subregular at $+\infty$. The result now follows from Lemma 5.1 and
Lemma 2.3 applied to appropriately chosen function
$h\colon\R^{n_1k_1\cdot\dots\cdot n_lk_l}\to \R_+$. \bbbox  \bs

\beginsection Section 6. Minmax hypercontractivity of norms of stable
random vectors.

In this section we apply the results in early sections to certain questions
concerning Gaussian and symmetric stable measures. 

The following lemma is a consequence of Kanter's inequality, (cf. Ledoux
and Talagrand (1991), p. 153)
which can be viewed as a concentration result similar to Levy's
inequalities. The formulation of the lemma below for Gaussian measures was
suggested by X. Fernique.

\ms
\proclaim Lemma 6.1 (Corollary of Kanter's inequality). Let $\nu$ be a
symmetric $\alpha$ stable measure with $0 < \alpha \le 2$ on a separable
Banach space $F$.
Then, for any $\kappa\ge 0$, any symmetric, convex set $B$ and any $y\in
F$, we have
$$
\nu(\kappa B+y)\le {\dsp {3\over 2} {\kappa^{\alf/2}\over \sqrt{1-\nu(B)}}}. $$

\n {\bf Proof.}
Let $\{X, X_i\}_i$ be i.i.d. symmetric $\alpha$ stable random variables
with $0 < \alpha \le 2$.
Take $N=[{\dsp { \kappa^{-\alpha}}}]$. Then using $N\kappa^\alf\le 1$ and
$(N+1)\kappa^\alf >1$, we have by Kanter's inequality $$ \eqalign{P(X-y
&\in \kappa B)
=P(\sum_{i=1}^N X_i - N^{1/\alf}y\in N^{1/\alf}\kappa B)\cr \le &{\dsp
{3\over 2} \Big({1\over
1+NP(X\notin N^{1/\alf} \kappa B)}\Big)^{1/2}}\le {\dsp {3\over
2}{\kappa^{\alf/2} \over P(X\notin B)^{1/2}}}\cr} $$ since $P(X \notin
N^{1/\alf} \kappa B)\ge P(X\notin B)$ and $(1+NP(X \notin B) )^{-1} \le
\kappa^\alf P^{-1}(X \notin B)$. This finishes the proof. \bbbox

\proclaim Lemma 6.2.
Let $\nu$ be a symmetric $\alpha$ stable measure with $0 < \alpha \le 2$ on
a separable, Banach space $F$. Then for any closed, symmetric, convex set
$B\subseteq F$, $y\in F$ and $\kappa\le 1$, $$ \nu(\kappa B+y)\le
R\kappa^{\alf/2}\nu(2B+y), $$ where $R=(3/ 2)
(\nu(B))^{-1}(1-\nu(B))^{-1/2}$.

\n {\bf Proof.} First consider $y\in B$. Then $\nu(B)\le\nu(2B+y)$ since
$B\subseteq 2B+y$. Thus, to conclude this case, one applies Lemma 6.1.

If $y\notin B$, then let $r=[{\dsp {\kappa^{-1}}-{2^{-1}}}]$. For
$k=0,1,\cdots,r$ the balls $\{y_k+\kappa B\}$ are disjoint and contained in
$y+2B$, where $y_k=(1-{\dsp {2\kappa \|y\|^{-1}}k})y$. By Anderson's
Theorem, it follows that
$$
\nu(y_k+\kappa B)\ge\nu(y+\kappa B)
$$
for $k=0,\cdots,r$.
Therefore, $\nu(\kappa B+y)\le{\dsp {(
r+1)^{-1}}}\nu(2B+y)\le\kappa\nu(2B+y)$. This proves the lemma, since $2\le
R$. \bbbox

\ms
\proclaim Proposition 6.3. Under the set up of Lemma 6.2, we have for each
$\kappa, t\le 1$,
$$\nu(\kappa tB)\le R '\kappa^{\alf/2}\nu(tB),$$ where $R '=3
(\nu(B/2))^{-1}(1-\nu(B/2))^{-1/2}$.

\n {\bf Proof.} \n Now for any $0\le t\le 1$, define the probability
measure $\nu_t$ by $\nu_t(C)=\nu(tC)=P(X/t \in C)$ where $X$ is the
symmetric $\alpha$ stable random variable with law $\nu$. Then $$
\nu*\nu_s(C)=P(X+X'/s \in C)=P((1+s^{-\alf})^{1/\alf}X \in C)=\nu_t(C), $$
where $t^{-\alf}=1+ s^{-\alf}$ and $X'$ is an independent copy of $X$.
Hence, by Lemma 6.2
$$
\eqalign{\nu(\kappa tB)&=\nu*\nu_s(\kappa B) =\int_F\nu(2\kappa B/2+y)\nu_s(dy)
\le (2\kappa)^{\alpha/2}R\int_F\nu(B+y)\nu_s(dy)\cr&\le
R'\kappa^{\alpha/2}\nu(tB).\cr} $$
\ms
\proclaim Theorem 6.4. Under the set up of Lemma 6.2, for each $b<1$, there
exists $R(b)$ such that for all $0\le t\le 1$, $$ \nu(tB)\le
R(b)t^{\alf/2}\nu(B), \hbox{ whenever } \nu(B)\le b. \leqno{(6.1)} $$

\n {\bf Proof.} Fix $B$ with $\nu(B)\le b$. Choose $s\ge 1$ so that
$\nu(sB)=b$. Now, apply the Proposition 6.3 with $\kappa=t$, to get
$$\nu(tB)=\nu(t \cdot {1 \over 2s}(2sB))\le R(b)t^{\alf/2}\nu({\dsp {1\over
2s}}(2sB))=R(b)t^{\alf/2}\nu(B), $$
where $R(b)=3 b^{-1}(1-b)^{-1/2}$.
\bbbox

\n {\bf Remark 6.5.}
In the case of $\alpha=2$ Theorem 6.4 was formulated in Szarek (1991),
Lemma 2.6,
where a weaker result,
which was sufficient for the main results of the paper, was actually proved. 
Recently,  Lata{\l}a proved that in the case of $\alf=2$, the conclusion 
of Theorem 6.4 holds whenever the measure $\nu$ is log concave.

\bigskip\n Related results on $\alpha$-stable measures can be found in
Lewandowski, Ryznar and \.Zak (1992). The key difference is that we need
the right hand side of (6.1) to involve $\mu(B)$ for all $B$ such that
$\mu(B) \le b$ and
the constant $R$ depending only on the number $b$.

\bigskip

\proclaim Corollary 6.6.
Let $0<\alpha\leq 2$, $0<p,q$. If $\alpha=2$ we assume that $q<\alpha$. If
$X$ is a $\alpha$-stable, symmetric vector in a Banach space then $\|X\|\in
\min\max {\cal H}_{p,q}(C)$ for some constant $C$ which depends only on
$\alpha$, $p$ and $q$.

\n {\bf Proof.} By the result of Szulga (1990) $\|X\|\in \max {\cal
H}_{p,q}(C_1)$ for some constant $C_1$ which depends only on $\alpha,p,q$.
Theorems 6.4 and 3.4 imply that $\|X\|\in \min {\cal H}_{p,q}(C_2)$ where
$C_2$ depends only on $\alpha,p,q$. Therefore $\|X\|\in \min\max {\cal
H}_{p,q}(C)$
with $ C= C_1\vee C_2$\bbbox

\proclaim Corollary 6.7.
Let $0<\alpha\leq 2$, $0<p,q$. If $\alpha=2$ we assume that $1<q<\alpha$.
Let $X_1,X_2,..,X_n$ be symmetric $\alpha$-stable, independent random
vectors in a Banach space. Let $h:\R^n_+ \rightarrow \R_+$ be a function as
in Theorem 5.3 which is $\lambda$-homogeneous for some $\lambda$. Then
$$h(\|X_1\|,\|X_2\|,..,\|X_n\|)\in \min\max {\cal H}_{p,q}(C)$$ and the
constant $C$ depends only on $\alpha,p,q$.

\n {\bf Proof.} By Corollary 6.6 a constant $\sigma$ can be found, which
depends only on $\alpha, p, q$ and such that the conclusion of Lemma 5.1
holds true for
$X=\|X_i\|$ for $i=1,2,..,n$. Now we can proceed as in the proof of Theorem
5.3.\bbbox

\proclaim Theorem 6.8.
Let $Y=\max_{l\le L} \|G\|_l$ and $X=\max_{l\le L} \|G_l\|_l$, where the
norms $\|\cdot\|_l$ were defined in
Section 1. If
$$
\|m_n(Y)\|_q \le C\|m_n(X)\|_q ,
\leqno(6.2)
$$
then for all $t \ge 0$
$$
P(Y \le ct) \ge P(X \le t)
$$
where the constant $c$ is independent of dimension $n$, the number $L$ and
the norms $\|\cdot\|_l$.

\n {\bf Proof.} In what follows, the statement that a constant is
independent of everything means that the constant is independent of
dimension $n$, the number $L$ and the norms $\|\cdot\|_l$, but can depend
on $p$, $q$ and other absolute constants. In order to apply Theorem 5.4, we
first need to show $$ \hbox{max hypercontractivity} \quad \|M_n(X)\|_q \le
C_1 \|M_n(X)\|_p, \leqno(6.3)
$$
$$
\hbox{min hypercontractivity} \quad \|m_n(X)\|_q \le C_2 \|m_n(X)\|_p,
\leqno(6.4)
$$
and
$$
\|M_n(Y)\|_q \le C_3 \|M_n(X)\|_q
\leqno(6.5)
$$
where constants $C_1$, $C_2$ and $C_3$ are independent of everything. In
particular, the constant $c$ is independent of everything.

\bigskip
To prove (6.3), note that $M_n(X)$ is a norm of Gaussian vectors. Thus by
the hypercontractivity of norms of Gaussian vectors (cf. for example,
Ledoux and Talagrand (1990), p. 60)
we obtain that $C_1$ can be a constant independent of everything.

\bigskip
To prove (6.4), we use Proposition 6.3 to check the condition (iii) of
Theorem 3.4. Taking $b=1/2$, by Lemma 2.1 and (6.3), we can take $\lmd$
close to $1$, but independent of everything (depending on $C_1$, $p$ and
$q$ only), such that $$ P(X \le t) \le b=1/2
$$
for all $t \le t_0= \lmd \| \max_{l\le L} \|G_l\|_l \|_p$. Note that $t_0$
here is not independent of everything, but we are interested in the
constant $C_2$ in (6.4) or $C$ in (i) of Theorem 3.4. By Proposition 6.3,
we have
$$
P(X\le st) \le RsP(X \le t) \quad \hbox{for all} \quad 0 <s <1 \quad
\hbox{and for all} \quad t \le t_0.
$$
For each $\vp >0$, taking $\tau=s=(\vp /R) \wedge (1/2)$, we obtain $$
P(X\le \tau t) \le \vp P(X\le t) \quad \hbox{\rm for \ all}\quad t\le t_0
$$ which implies (6.4) by Theorem 3.4 with $C_2$ independent of everything.
To see that the constant $C_2$ is
independent of everything, we only need to follow the part of the proof
from (iii) $\Rightarrow$ (iv) and (iv) $\Rightarrow$ (i) in Theorem 3.4. It
is clear that $\tau$ and $r$ in the proof is independent of everything and $$
C=\sigma^{-1}
=\max\big\{{ \|X\|_q / (E((\tau t_0)\wedge X)^p)^{1/p}}, (\tau
r)^{-1}\big\}. $$ Now note that
$$
\eqalign{
(E((\tau t_0)\wedge Y)^p)^{1/p}
&\ge (\tau t_0)P^{1/p}(X \ge \tau t_0) \cr &=\tau \lmd \|X\|_p P^{1/p} (X
\ge \lmd \tau \| X\|_p) \cr &\ge \tau \lmd \|X\|_p \cdot P^{1/p} (X \ge
\lmd \tau \| X\|_p) } $$ and by Lemma 2.1 and (6.3)
$$
P (X \ge \lmd \tau \| X\|_p)
\ge ((1-(\lmd\tau)^p)C_1^{-p})^{q/(q-p)}. $$ Thus
$$
\eqalign{
C&=\max\big\{{ \|X\|_q / (E((\tau t_0)\wedge X)^p)^{1/p}}, (\tau
r)^{-1}\big\}\cr &\ge \max\big\{ \|X\|_q / (\|X\|_p \tau\lmd (1-(\lmd\tau)^p
)C_1^{-p})^{q/p(q-p)}) , (\tau r)^{-1}\big\} \cr
&\ge \max\big\{ C_1 (\tau\lmd (1-(\lmd\tau)^p)C_1^{-p})^{q/p(q-p)})^{-1} ,
(\tau r)^{-1}\big\}\cr
&=C_2 }
$$
and it is clear that $C_2$ is independent of everything.

Finally, (6.5) follows from Slepian's lemma, see (1.2) and remarks following it.

Now we can apply Theorem 5.4 with (6.3), (6.4) and (6.5) in hand. We only
need to show that
$\sigma$ in Lemma 5.1 is independent of everything in view of the proof of
Theorem 5.4. To check that $$
\sigma=\min\big\{{\alpha t_0 /\|Y\|_q},
{\alpha/ D}, {r/ D}, r\tau \big\}
$$
in Lemma 5.1, is similar to the proof of (6.4) and is omitted here. Thus we
finished the proof of Theorem 6.8. \bbbox \bigskip

As a consequence of Theorem 6.8, we have the following modified correlation
inequality for centered Gaussian measure.

\proclaim Corollary 6.9.
(modified correlation inequality)
Assume (6.2) holds. Then there exists an absolute constant $\alf$ such that
$$\mu\bigl(\alf(\bigcap_{i=1}^l A_i)\bigr)\ge\prod_{i=1}^l\mu(A_l) \leqno(6.6) $$
for any centered Gaussian measure $\mu$ and any convex, symmetric sets
$A_i$, $1 \le i \le l$ and $l \ge 1$.

\bigskip
\n {\bf Remark 6.10.}
Note that the original correlation conjecture (1.1) implies (6.6) with
constant $c=1$.
\bs

\beginsection Section 7. Final remarks and some open problems.

In this section we mention a few results and open problems that are closely
related to the main results in this paper. At first, we give a very simple
proof of the following result.

\proclaim Proposition 7.1.
For $0<p<q<\infty$, if there exists a constant $C$, such that for all $n$,
$$ \|M_n(X)\|_q\le C\|M_n(X)\|_p
\leqno(7.1)
$$
then the following are equivalent:
\item{(i)} There exists a constant $C$, such that for all $n$, $$
\|M_n(Y)\|_p\le C\|M_n(X)\|_p ;
\leqno(7.2)
$$
\item{(ii)} There exists a constant $C$ such that $$ P(Y >t) \le C \cdot
P(X>t). \leqno(7.3)
$$

\n {\bf Proof.}
It follows from de~la~Pe\~na, Montgomery-Smith and Szulga (1994) that the
hypercontractivity of $X$, (7.1), and the domination relation (7.2) imply
the tail domination (7.3).
So we only need to show that (ii) implies (i). Without loss of generality,
we assume $C>1$.
Let $\delta$ be an independent random variable with $$ P(\delta =1)=1/C,
\quad P(\delta =0)=1-1/C. $$ Then for all $n$ and all $t\ge 0$
$$
\eqalign{ P( M_n(\delta Y) < t )
= & P^n( \delta Y < t )
=(1-P( \delta Y \ge t ))^n
=(1-C^{-1}P( Y \ge t ))^n \cr
& \ge (1-P( X \ge t ))^n
=P( M_n(X) < t ) }
$$
which implies $\|M_n(\delta Y)\|_p\le \|M_n(X)\|_p $. On the other hand, we
have $$
E_Y E_\delta \max_{1 \le i \le n} ( \delta_i Y_i^p) \ge E_Y \max_{1 \le i
\le n} E_\delta ( \delta_i Y_i^p) =
C^{-1}EY_i^p.
$$
which finishes the proof.\bbbox
\bs

Our next proposition is related to Theorem 6.4. \proclaim Proposition 7.2.
Let $\mu$ be a probability measure on a separable, Banach space $F$. Then
for any closed, symmetric, convex set $B\subseteq F$, the following are
equivalent:
\item{(i)} For each $b<1$, there exist $R=R(b)>1$ and $\beta>0$ such that
for all $0\le t\le 1$,
$$
\mu(tB)\le R(b)t^\beta\mu(B), \hbox{ whenever } \mu(B)\le b. \leqno{(7.4)} $$
\item{(ii)} For each $b<1$, there exists $r=r(b)<1$ such that for all $t
\ge 0$, $$
\int_0^t \mu(sB)ds \le r(b)t\mu(tB), \hbox{ whenever } \mu(B)\le b.
\leqno{(7.5)}
$$

\n {\bf Proof.}
We first prove (i) implies (ii). Take $\delta<1$ such that
$\delta^\beta=(1+\beta)^{-1} R^{-1}$.
Then we have
$$
\eqalign{ \int_0^t \mu(sB)ds
& =\int_0^{\delta t}\mu(sB)ds +\int_{\delta t}^t\mu(sB)ds \cr & \le
\int_0^{t} \delta \mu(\delta (sB))ds +(1-\delta )t\mu(tB) \cr & \le R
\delta^{1+\beta} \int_0^{t} \mu(sB)ds +(1-\delta)t\mu(tB) \cr & \le
(1-\delta+R\delta^{1+\delta})t\mu(tB) \cr & = (1-\beta (1+\beta)^{-1}
\delta)t\mu(tB) .} $$
Thus (ii) holds with $r(b)=1-\beta (1+\beta)^{-1} \delta <1$.

\n To prove (ii) implies (i), note that for $\delta=1-r^{1/2} >0$ $$
\int_0^t \mu(sB)ds \ge \int_{\delta t}^t \mu(sB)ds \ge (1-\delta)t
\mu(\delta tB) =r^{1/2}t\mu(\delta tB) $$ and thus from (7.5), for any
$t>0$,
$$
\mu(\delta tB) \le r^{1/2} \mu(tB).
\leqno{(7.6)}
$$
Now for any $0<s<1$, pick $k\ge 0$ such that $\delta^{k+1} \le s
<\delta^k$. We have by the iteration of (7.6) that $$
\mu(sB) \le \mu(\delta^k B) \le r^{1/2} \mu(\delta^{k-1} B) \le r^{k/2}
\mu( B). $$
Thus by using $k+1 \ge \log s/\log \delta$, we obtain $$ \mu(sB)\le r^{k/2}
\mu( B) \le r^{-1/2} s^{\log r/(2\log \delta)} \mu( B) $$ which finishes
the proof with $R=r^{-1/2}>1$ and $\beta=\log r/(2\log \delta)>0$.
\bbbox
\bs

\n There are many questions related to this work. Let us only mention a few
here.

\proclaim Conjecture 7.3.
The best min-hypercontractive constant in (6.4) with $Y=\|X\|$ for
symmetric Gaussian vectors $X$ in any Banach space is $$ C ={\Gamma^{1/q}
(q) \over \Gamma^{1/p} (p)} . $$ \bs

The constant follows from the small ball estimates, $P(|X|<s)\sim K \cdot
s$ as $s \to 0$, of one-dimensional Gaussian random variable $X$. Note that
if $\beta >1$ and $P(|X|<s)\sim K \cdot s^\beta$ as $s \to 0$, then the
resulting constant in this case is smaller. Thus the conjecture looks
reasonable in view of Proposition 6.3.

A related question is that under hyper-max condition, what can we say about
a nontrivial lower bound for $\|M_{k+1}\|_p / \|M_k\|_p $, and in
particular, in the Gaussian case. This maybe useful in proving the
conjecture.

A result of Gordon (1987) compares the expected minima of maxima for, in
particular, Gaussian processes. We mention this here because a version of
Gordon's results could perhaps be used to prove the next Conjecture. Note
that if the conjecture holds, then the modified correlation inequality $
C_\alf$ holds.

\proclaim Conjecture 7.4.
Let $G$, $G_l$ and norm $\| \cdot \|_l$ be as in Section 1. If
$Y=\max_{l\le L} \|G\|_l$ and
$X=\max_{l\le L} \|G_l\|_l$, then $$
\|m_n(Y)\|_q \le C\|m_n(X)\|_q .
$$
\bs

Our next conjecture is related to stable measures. It is a stronger
statement than our Proposition 6.4 and holds for the symmetric Gaussian
measures.

\proclaim Conjecture 7.5.
Let $\nu$ be a symmetric $\alpha$ stable measure with $0 < \alpha \le 2$ on
a separable, Banach space $F$. Then for any closed, symmetric, convex set
$B\subseteq F$ and for each $b<1$,
there exists $R(b)$ such that for all $0\le t\le 1$, $$ \nu(tB)\le
R(b)t\nu(B), \hbox{ whenever } \nu(B)\le b. $$ \bs

Note also that the following S-conjecture will provide the best constant
$R(b)$ for the symmetric Gaussian measures in Theorem 6.4, (6.1). See
Kwapie\'{n} and Sawa (1993) for history and a proof of the conjecture for
1-unconditional set. \bs

\proclaim S-Conjecture.
Let $\mu$ be a symmetric Gaussian measure on a separable Banach space $F$.
Then for any closed, symmetric, convex set $B\subseteq F$, $$ \mu(\lambda
B) \ge \mu(\lambda S)
$$
for each $\lambda >1$ and each symmetric slab $S$ in $F$ such that $\mu(S)
\ge \mu(B)$.

\vfill\eject

\newcount\probno
\def\r{\par\global\advance\seriesanswer by 1
\par\noindentitem{(\romannumeral\seriesanswer)}} \newcount\seriesanswer
\def\ms{\medskip}
\centerline{\bf Bibliography}
\bigskip

\par\noindent
Anderson, T.W. (1955). The integral of a symmetric unimodal function over a
symmetric convex set and some probability inequalities. {\it Proc. Amer.
Math. Soc.} {\bf 6}, 170--176. \ms

\noindent
Asmar, N.~H., Montgomery-Smith, S.J. (1994). On the distribution of Sidon
series. {\it Ark. Mat.} {\bf 31}, 13--26. \ms
\par\noindent
Bingham, N.~H., Goldie, C.~M. and Teugels, J.~L. (1987). Regular variation. {\it Cambridge University Press,} 491 pp.
\ms\par\noindent
De~la~Pe\~na, V., Montgomery-Smith, S.J. and Szulga, J. (1994). Contraction
and decoupling inequalities for multilinear forms and U-statistics. {\it
Ann. of Prob.} {\bf 22}, 1745--1765. \ms
\par\noindent
Gordon, Y. (1987).
Elliptically contoured distributions. {\it Prob. Theory Rel. Fields} {\bf
76}, 429--438. \ms
\par\noindent
Khatri, C.G. (1967). On certain inequalities for normal distributions and
their applications to simultaneous confidence bounds. {\it Ann. Math.
Statist.} {\bf 38}, 1853--1867. \ms

\par\noindent
Kuelbs, J., Li, W.V. and Shao, Q. (1995). Small ball estimates for
fractional Brownian motion under H\"older norm and Chung's functional LIL.
{\it J. Theor. Prob.} {\bf 8}, 361-386. \ms

\par\noindent
Kwapie\'{n}, S. and Sawa, J. (1993),
On some conjecture concerning Gaussian measures of dilatations of convex
symmetric sets,
{\it Studia Math.} {\bf 105}, 173--187. \ms

\par\noindent
Ledoux, M. and Talagrand, M. (1991).
{\it Probability on Banach Spaces}, Springer, Berlin. 
\ms \par\noindent
Lewandowski, M., Ryznar, M. and \.Zak, T. (1992). Stable measure of a small
ball,
{\it Proc. Amer. Math. Soc.}, {\bf 115}, 489--494. \ms
\par\noindent
Li, W.V. and Shao, Q. (1995).
The existence of the small ball constant for sup-norm of fractional
Brownian motion under the correlation conjecture. (notes). \ms
\par\noindent
Marcus, M.~B. and Shepp, L. (1972). Sample behavior of Gaussian processes.
{\it Proc. of the Sixth Berkeley Symposium on Math. Statist. and Prob.}
{\bf vol. 2}, 423--441 \ms
\par\noindent
Pitt, L.D. (1977).
A Gaussian correlation inequality for symmetric convex sets. {\it Ann.
of Prob.} {\bf 5}, 470--474. \ms
\par\noindent
Rychlik, E. (1993). Some necessary and sufficient conditions of $(p,q)$
hypercontractivity. Preprint. \ms
\par\noindent
Schechtman, G., Schlumprecht, T. and Zinn, J. (1995). On the Gaussian
measure of the intersection of symmetric, convex sets. Preprint. \ms

\par\noindent
S\v id\' ak, Z. (1967) Rectangular confidence regions for the means of
multivariate normal distributions. {\it J. Amer. Statist. Assoc.} {\bf 62},
626--633. \ms
\par\noindent
S\v id\' ak, Z. (1968). On multivariate normal probabilities of rectangles:
their dependence on correlations. {\it Ann. Math. Statist.} {\bf 39},
1425--1434. \ms

\par\noindent
Slepian, D. (1962). The one-sided barrier problem for Gaussian noise. {\it
Bell System Tech. J. 41}, 463--501. \ms

\par\noindent
Szarek, S. (1991).
Conditional numbers of random matrices.
{\it J. Complexity}, {\bf 7}, 131--149.
\ms

\par\noindent
Szarek, S. and Werner, E. (1995). Personal Communication. \ms

\par\noindent
Szulga, J. (1990). A note on hypercontractivity of stable random variables.
{\it Ann. of Prob.} {\bf 18}, 1746--1758. \ms

\par\noindent
Talagrand, M. (1994).
The small ball problem for the Brownian sheet.  {\it Ann. of Probab.} {\bf
22}, 1331--1354.
\ms

\par\noindent
Tong, Y. L. (1980).
{\it Probability Inequalities in Multivariate Distributions}, Academic
Press, New York.
\vfill\eject
\settabs
\+ College Station, TX 77843--3368, USAAAAA& Department of Theoretical
Mathematicsasssssssss\cr \+P. Hitczenko&S. Kwapie\'{n} \cr \+Department of
Mathematics,&Institute of Mathematics\cr \+North Carolina State
University&Warsaw University
\cr \+Raleigh, NC 27695--8205, USA&Banacha 2, 02-097 Warsaw, Poland\cr
\medskip \settabs \+ College Station, TX 77843--3368, USAAAAA&Department of
Theoretical
Mathematicsasssssssss\cr \+W. V. Li&G. Schechtman\cr \+Department of
Mathematics,&Department of Theoretical Mathematics \cr \+University of
Delaware&The Weizmann Institute of Science \cr \+Newark, DE 19716,
USA&Rehovot, Israel\cr
\medskip\settabs \+ College Station, TX 77843--3368, USAAAAA& Department of
Theoretical Mathematicsasssssssss\cr
\+T. Schlumprecht& J. Zinn\cr
\+ Department of Mathematics& Department of Mathematics\cr \+ Texas A{\&M}
University&Texas A{\&M} University\cr \+ College Station, TX 77843--3368,
USA&College Station, TX 77843--3368, USA \cr
\bye